\documentclass[a4paper,12pt]{article}
\usepackage[cp1251]{inputenc}
\usepackage[russian]{babel}
\usepackage{amsfonts, amssymb, amsmath, amsthm, amscd}
\usepackage{cite}
\author{A.A. Vasil'eva}
\title{Estimates for norms of two-weighted summation operators on trees for $1<p<q<\infty$}
\date{}
\begin{document}

\maketitle

\newenvironment{Biblio}{%
                  \renewcommand{\refname}{\footnotesize REFERENCES}%
                  \begin{thebibliography}{99}%
                  \itemsep -0.2em}%
                  {\end{thebibliography}}

\def\inff{\mathop{\smash\inf\vphantom\sup}}
\renewcommand{\le}{\leqslant}
\renewcommand{\ge}{\geqslant}
\newcommand{\sgn}{\mathrm {sgn}\,}
\newcommand{\inter}{\mathrm {int}\,}
\newcommand{\dist}{\mathrm {dist}}
\newcommand{\supp}{\mathrm {supp}\,}
\newcommand{\R}{\mathbb{R}}
\renewcommand{\C}{\mathbb{C}}
\newcommand{\Z}{\mathbb{Z}}
\newcommand{\N}{\mathbb{N}}
\newcommand{\Q}{\mathbb{Q}}
\theoremstyle{plain}
\newtheorem{Trm}{Theorem}
\newtheorem{trma}{Theorem}
\newtheorem{Def}{Definition}
\newtheorem{Cor}{Corollary}
\newtheorem{Lem}{Lemma}
\newtheorem{Rem}{Remark}
\newtheorem{Sta}{Proposition}
\newtheorem{Exa}{Example}
\renewcommand{\proofname}{\bf Proof}
\renewcommand{\thetrma}{\Alph{trma}}
\renewcommand{\abstractname}{Abstract}

\begin{abstract}
In this paper, estimates for norms of weighted summation operators
(discrete Hardy-type operators) on a tree are obtained for
$1<p<q<\infty$ and for arbitrary weights and trees.
\end{abstract}

\section{Introduction}

First we give some notation.

Throughout this paper we consider graphs ${\cal G}$ with finite or
countable vertex set, which will be denoted by ${\bf V}({\cal
G})$. Also we suppose that the graphs have neither multiple edges
nor loops. Given a function $f:{\bf V}({\cal G})\rightarrow \R$
and a number $1<p<\infty$, we set
$$
\|f\|_{l_p({\cal G})}=\left ( \sum \limits _{\xi \in {\bf V}
({\cal G})}|f(\xi)|^p\right )^{1/p}.
$$
Denote by $l_p({\cal G})$ the space of functions $f:{\bf V}({\cal
G})\rightarrow \R$ with finite norm $\|f\|_{l_p({\cal G})}$.

Let ${\cal T}=({\cal T}, \, \xi_0)$ be a tree rooted at $\xi_0$.
We introduce a partial order on ${\bf V}({\cal T})$ as follows: we
say that $\xi'>\xi$ if there exists a simple path $(\xi_0, \,
\xi_1, \, \dots , \, \xi_n, \, \xi')$ such that $\xi=\xi_k$ for
some $k\in \{0, \,\dots, \,  n\}$; by the distance between $\xi$
and $\xi'$ we mean the quantity $\rho_{{\cal T}}(\xi, \,
\xi')=\rho_{{\cal T}}(\xi', \, \xi) =n+1-k$. In addition, we set
$\rho_{{\cal T}}(\xi, \, \xi)=0$. If $\xi'>\xi$ or $\xi'=\xi$, we
write $\xi'\ge \xi$. For $j\in \Z_+$ and $\xi\in {\bf V}({\cal
T})$, let
$$
\label{v1v} {\bf V}_j ^{{\cal T}}(\xi):= \{\xi'\ge\xi:\;
\rho_{{\cal T}}(\xi, \, \xi')=j\}.
$$
Given $\xi\in {\bf V}({\cal T})$, we denote by ${\cal
T}_\xi=({\cal T}_\xi, \, \xi)$ the subtree in ${\cal T}$ with the
vertex set
$$
{\bf V}({\cal T}_\xi)=\{\xi'\in {\bf V}({\cal T}):\xi'\ge \xi\}.
$$

Let ${\mathbf W}\subset {\mathbf V}({\mathcal T})$. We say that
${\mathcal G}\subset {\mathcal T}$ is a maximal subgraph on the
vertex set ${\mathbf W}$ if ${\mathbf V}({\mathcal G})={\mathbf
W}$ and if any two vertices $\xi'$, $\xi''\in {\mathbf W}$ that
are adjacent in ${\mathcal T}$ are also adjacent in ${\mathcal
G}$. Given $\xi$, $\xi'\in {\bf V}({\cal T})$, $\xi\le \xi'$, we
denote by $[\xi, \, \xi']$ the maximal subgraph on the vertex set
$\{\eta\in {\bf V}({\cal T}):\; \xi\le \eta\le \xi'\}$.

Let ${\cal G}$ be a subgraph in $({\cal T}, \, \xi_0)$. Denote by
${\bf V}_{\max} ({\cal G})$ and ${\bf V}_{\min}({\cal G})$ the set
of maximal and minimal vertices in ${\cal G}$, respectively.

Let $({\cal T}, \, \xi_0)$ be a tree, and let $u$, $w:{\bf
V}({\cal T})\rightarrow [0, \, \infty)$ be weight functions.
Define the summation operator $S_{u,w,{\cal T}}$ by
$$
S_{u,w,{\cal T}}f(\xi) = w(\xi)\sum \limits _{\xi'\le \xi}
u(\xi')f(\xi'), \quad \xi \in {\bf V}({\cal T}), \quad f:{\bf
V}({\cal T}) \rightarrow \R.
$$
Given $1< p, \, q< \infty$, by $\mathfrak{S}^{p,q}_{{\cal T},u,w}$
we denote the operator norm of $S_{u,w,{\cal T}}:l_p({\cal T})
\rightarrow l_q({\cal T})$, which is the minimal constant $C$ in
the inequality
$$
\left(\sum \limits_{\xi \in {\bf V}({\cal T})} w^q(\xi) \left|
\sum \limits _{\xi'\le \xi}u(\xi')f(\xi')\right|^q\right)^{1/q}
\le C\left(\sum \limits_{\xi \in {\bf V}({\cal T})}
|f(\xi)|^p\right)^{1/p}, \;\; f:{\bf V}({\cal T})\rightarrow \R.
$$

Let $X$, $Y$ be arbitrary sets, $f_1$, $f_2:\ X\times Y\rightarrow
\R_+$. We write $f_1(x, \, y)\underset{y}{\lesssim} f_2(x, \, y)$
(or $f_2(x, \, y)\underset{y}{\gtrsim} f_1(x, \, y)$) if, for any
$y\in Y$, there exists $c(y)>0$ such that $f_1(x, \, y)\le
c(y)f_2(x, \, y)$ for each $x\in X$; $f_1(x, \,
y)\underset{y}{\asymp} f_2(x, \, y)$ if $f_1(x, \, y)
\underset{y}{\lesssim} f_2(x, \, y)$ and $f_2(x, \,
y)\underset{y}{\lesssim} f_1(x, \, y)$.

\begin{Trm}
\label{sum_oper} Let $({\cal A}, \, \xi_0)$ be a tree, and let
$u$, $w:{\bf V}({\cal A}) \rightarrow \R_+$. Suppose that
$1<p<q<\infty$. Then
$$
\mathfrak{S}^{p,q}_{{\cal A},u,w} \underset{p,q}{\asymp} \sup
_{\xi\in {\bf V}({\cal A})} \|u\|_{l_{p'}([\xi_0, \, \xi])}
\|w\|_{l_q({\cal A}_\xi)}.
$$
\end{Trm}

In \cite{vas_har_tree} this result was proved under some
restrictions on weights (see Theorems 1.2 and 3.6).

\smallskip

Given $f:{\bf V}({\cal A}) \rightarrow \R$, we set
$$
\|f\|_{l_q(l_p({\cal A}))} = \left(\sum \limits _{j=0}^\infty
\left(\sum \limits _{\xi \in {\bf V}_j^{{\cal A}}(\xi_0)}
|f(\xi)|^p\right) ^{\frac qp}\right)^{\frac 1q}.
$$
By $\hat{\mathfrak{S}}^{p,q}_{{\cal A},u,w}$ we denote the
operator norm of $S_{u,w,{\cal A}}: l_q(l_p({\cal A})) \rightarrow
l_q({\cal A})$. If $p\le q$, then $\|f\|_{l_p({\cal A})}\ge
\|f\|_{l_q(l_p({\cal A}))}$ and
\begin{align}
\label{spqh}\hat{\mathfrak{S}}^{p,q}_{{\cal A},u,w} \ge
\mathfrak{S} ^{p,q}_{{\cal A},u,w}.
\end{align}

Let
\begin{align}
\label{ujwj} u(\xi)=u_j, \quad w(\xi) =w_j, \quad \xi \in {\bf
V}^{\cal A}_j(\xi_0).
\end{align}
In addition, we suppose that there exist a number $C_*\ge 1$ and a
function $S:\Z_+ \rightarrow (0, \, \infty)$ satisfying the
following conditions:
\begin{align}
\label{cvj} C_*^{-1} \frac{S(j')}{S(j)}\le {\rm card}\, {\bf
V}_{j-j'}^{{\cal A}}(\xi)\le C_*\frac{S(j')}{S(j)}, \quad \xi \in
{\bf V}^{\cal A}_j(\xi_0), \quad j'\ge j;
\end{align}
\begin{align}
\label{s01} S(0)=1;
\end{align}
there exist $R_0\ge R>1$ such that
\begin{align}
\label{sj} R_0\ge \frac{S(j+1)}{S(j)} \ge R.
\end{align}

\begin{Trm}
\label{lpq} Suppose that $1<p<q<\infty$ and conditions
(\ref{ujwj}), (\ref{cvj}), (\ref{s01}), (\ref{sj}) hold. Then
$$\hat{\mathfrak{S}}^{p,q}_{{\cal A},u,w} \underset{p,q, C_*,R,
R_0}{\asymp} \sup _{j\in \Z_+} u_j\left(\sum \limits _{i\ge j}
w_i^q \frac{S(i)}{S(j)}\right)^{\frac 1q} \underset{p,q, C_*,R,
R_0}{\asymp} \sup _{\xi \in {\bf V}({\cal A})} u(\xi)
\|w\|_{l_q({\cal A}_\xi)}.$$
\end{Trm}

\section{Proof of Theorem \ref{sum_oper}}

The lower estimate for $\mathfrak{S}^{p,q}_{{\cal A},u,w}$ was
obtained in \cite[Lemma 3.3]{vas_har_tree}. In addition, the
following result was proved (see \cite[Lemma 3.1]{vas_har_tree}).

\begin{Lem}
\label{sigma} Let $1<p<q<\infty$. Then there exists
$\sigma=\sigma(p, \, q)\in \left(0, \, \frac 18\right)$ with the
following property: if $({\cal A}, \, \xi_0)$ is a tree with
finite vertex set, $u$, $w:{\bf V}({\cal A}) \rightarrow (0, \,
\infty)$,
$$
\frac{\|w\|_{l_q({\cal A}_\eta)}}{\|w\|_{l_q({\cal A}_\xi)}} \le
\sigma \quad \text{for any}\quad \xi \in {\bf V}({\cal A}), \quad
\eta \in {\bf V}_1^{\cal A}(\xi),
$$
then $\mathfrak{S}^{p,q}_{{\cal A},u,w} \underset{p,q}{\asymp}
\sup _{\xi\in {\bf V}({\cal A})} u(\xi) \|w\|_{l_q({\cal
A}_\xi)}$.
\end{Lem}
Lemma \ref{sigma} was proved by induction; here the discrete
analogue of Evans -- Harris -- Pick theorem \cite{ev_har_pick} was
applied.

Let $\{{\cal T}_j\}_{j\in \N}$ be a family of subtrees in ${\cal
T}$ such that ${\bf V}({\cal T}_j)\cap {\bf V}({\cal T}_{j'})
=\varnothing$ for $j\ne j'$ and $\cup _{j\in \N} {\bf V}({\cal
T}_j) ={\bf V}({\cal T})$. Then $\{{\cal T}_j\} _{j\in \N}$ is
called a partition of the tree ${\cal T}$. Let $\xi_j$ be the
minimal vertex of ${\cal T}_j$. We say that the tree ${\cal T}_s$
succeeds the tree ${\cal T}_j$ (or ${\cal T}_j$ precedes the tree
${\cal T}_s$) if $\xi_j<\xi_s$ and $$\{\xi\in {\cal T}:\; \xi_j\le
\xi<\xi_s\} \subset {\bf V}({\cal T}_j).$$

\renewcommand{\proofname}{\bf Proof of Theorem \ref{sum_oper}}
\begin{proof}
By B. Levi's theorem, without loss of generality we may assume
that the tree ${\cal A}$ has a finite vertex set. In addition, we
may consider only strictly positive weight functions $u$, $w$.

Let $\sigma=\sigma(p,\, q)\in \left(0, \, \frac 18\right)$ be as
defined in Lemma \ref{sigma}. Given $\xi \in {\bf V}({\cal A})$,
we denote
\begin{align}
\label{vxisig} {\bf V}_{\xi,\sigma} =\{\eta\ge \xi:\;
\|w\|_{l_q({\cal A}_\eta)} \ge \sigma \|w\|_{l_q({\cal A}_\xi)}\}.
\end{align}
Then $\xi\in {\bf V}_{\xi,\sigma}$. Denote by ${\cal
A}_{\xi,\sigma}$ the maximal subgraph in ${\cal A}$ on the vertex
set ${\bf V}_{\xi,\sigma}$. Notice that ${\cal A}_{\xi,\sigma}$ is
a tree.

We claim that
\begin{align}
\label{uw} \|u\|_{l_{p'}({\cal A}_{\xi,\sigma})} \|w\|_{l_q({\cal
A}_\xi)} \underset{p,q}{\lesssim} \sup _{\eta \in {\bf V}({\cal
A})} \|u\|_{l_{p'}([\xi_0, \, \eta])} \|w\|_{l_q({\cal A}_\eta)}.
\end{align}

Indeed, let $\eta \in {\bf V}({\cal A}_{\xi,\sigma})$. Then there
exists a vertex $\zeta_\eta \in {\bf V}_{\max}({\cal
A}_{\xi,\sigma})$ such that $\eta \in [\xi, \, \zeta_\eta]$.
Indeed, otherwise for any vertex $\zeta\in {\bf V}({\cal
A}_{\xi,\sigma})\cap {\bf V}({\cal A}_\eta)$ the set ${\bf
V}({\cal A}_{\xi,\sigma})\cap {\bf V}_1^{\cal A}(\zeta)$ is
nonempty. Hence, we can construct an infinite chain
$\zeta_1<\zeta_2<\zeta_3<\dots$, $\zeta_j\in {\bf V}({\cal
A}_{\xi,\sigma})$. This contradicts with the assumption that the
set ${\bf V}({\cal A})$ is finite.

Thus, ${\bf V}({\cal A}_{\xi,\sigma}) =\cup _{\zeta\in {\bf
V}_{\max}({\cal A}_{\xi,\sigma})} [\xi, \, \zeta]$. Therefore,
\begin{align}
\label{ulp} \|u\|_{l_{p'}({\cal A}_{\xi,\sigma})}^{p'} \le
\sum_{\zeta\in {\bf V}_{\max}({\cal A}_{\xi,\sigma})} \sum \limits
_{\xi\le \eta\le \zeta} u^{p'}(\eta).
\end{align}
Let us prove that
\begin{align}
\label{cvmax} {\rm card}\, {\bf V}_{\max}({\cal A}_{\xi,\sigma})
\underset{p,q}{\lesssim} 1.
\end{align}
Indeed, if $\zeta$, $\zeta'\in {\bf V}_{\max}({\cal
A}_{\xi,\sigma})$, then these vertices are incomparable and ${\bf
V}({\cal A}_{\zeta}) \cap {\bf V}({\cal A}_{\zeta'})
=\varnothing$. Hence,
$$
\|w\|^q_{l_q({\cal A}_\xi)} \ge \sum_{\zeta\in {\bf V}
_{\max}({\cal A}_{\xi,\sigma})} \|w\|^q_{l_q({\cal A}_\zeta)}
\stackrel{(\ref{vxisig})}{\ge} {\rm card}\, {\bf V}_{\max}({\cal
A}_{\xi,\sigma}) \sigma^q \|w\|^q_{l_q({\cal A}_\xi)};
$$
i.e., ${\rm card}\, {\bf V}_{\max}({\cal A}_{\xi,\sigma})\le
\sigma^{-q}$. From (\ref{ulp}) and (\ref{cvmax}) it follows that
there exists a vertex $\zeta_*\in {\bf V}_{\max}({\cal
A}_{\xi,\sigma})$ such that
$$
\|w\|_{l_q({\cal A}_\xi)}^{p'}\|u\|_{l_{p'}({\cal
A}_{\xi,\sigma})}^{p'} \underset{p,q}{\lesssim} \|w\|_{l_q({\cal
A}_\xi)}^{p'}\sum \limits _{\xi\le \eta\le \zeta_*}
u^{p'}(\eta)\stackrel{(\ref{vxisig})}{\le}
$$
$$
\le \sigma^{-p'}\|w\|^{p'}_{l_q({\cal A}_{\zeta_*})} \sum \limits
_{\xi\le \eta\le \zeta_*} u^{p'}(\eta).
$$
This implies (\ref{uw}).

Let us construct a partition of the tree ${\cal A}$ into subtrees
$({\cal A}_m, \, \hat\xi_m)$, $1\le m\le m_*$, such that
\begin{align}
\label{am} {\cal A}_m={\cal A}_{\hat \xi_m,\sigma}, \quad 1\le
m\le m_*.
\end{align}
To this end we construct a family of partitions of the tree ${\cal
A}=\tilde {\cal A}_k\sqcup {\cal G}_k$, $0\le k\le k_*$; here
$\tilde {\cal A}_k$ is a subtree in ${\cal A}$ rooted at $\xi_0$,
${\bf V}(\tilde {\cal A}_k) =\sqcup _{m=1}^{m_k} {\bf V}({\cal
A}_m)$, and ${\cal G}_k$ is a disjoint union of trees ${\cal
A}_{\eta_{k,j}}$, $1\le j\le j_k$. In addition, $m_{k+1}>m_k$.
\begin{enumerate}
\item We set ${\bf V}(\tilde{\cal A}_0)= \varnothing$, ${\cal
G}_0={\cal A}$, $m_0=0$.

\item Suppose that a partition ${\cal A}=\tilde {\cal A}_k\sqcup
{\cal G}_k$ is constructed and ${\bf V}({\cal G}_k) \ne
\varnothing$. Denote by $\tilde {\cal A}_{k+1}$ the maximal
subgraph on the vertex set ${\bf V}(\tilde{\cal A}_k) \cup
\left(\cup _{1\le j\le j_k} {\bf V} ({\cal A} _{\eta_{k,j},
\sigma})\right)$, and by ${\cal G}_{k+1}$, the maximal subgraph on
the vertex set ${\cal G}_k \backslash \left(\cup _{1\le j\le j_k}
{\cal A}_{\eta_{k,j},\sigma}\right)$. Let $\{{\cal A}_m\}_{1\le
m\le m_{k+1}} =\{{\cal A}_m\}_{1\le m\le m_{k}} \cup \{{\cal
A}_{\eta_{k,j},\sigma}\}_{1\le j\le j_k}$. If ${\bf V}({\cal
G}_{k+1})=\varnothing$, then we stop.
\end{enumerate}

From (\ref{vxisig}) and (\ref{am}) it follows that
\begin{align}
\label{wwsig} \frac{\|w\|_{l_q({\cal A}_{\hat \xi
_k})}}{\|w\|_{l_q({\cal A}_{\hat \xi_m})}}<\sigma \quad
\text{if}\quad {\cal A}_k \;\text{succeeds}\; {\cal A}_m.
\end{align}

We define the tree ${\cal D}$ with the vertex set $\{\hat
\xi_m\}_{1\le m\le m_*}$ as follows: we write $\hat \xi_l\in {\bf
V}_1^{{\cal D}}(\hat \xi_m)$ if and only if the tree ${\cal A}_l$
succeeds the tree ${\cal A}_m$. Let $\hat u$, $\hat w:{\bf
V}({\cal D}) \rightarrow (0, \, \infty)$,
\begin{align}
\label{hatuxim} \hat u(\hat \xi_m) =\|u\|_{l_{p'}({\cal A}_{\hat
\xi_m,\sigma})}, \quad \hat w(\hat \xi_m) =\|w\|_{l_{q}({\cal
A}_{\hat \xi_m,\sigma})}.
\end{align}
Then
\begin{align}
\label{spq} \mathfrak{S}^{p,q}_{{\cal A},u,w} \le
\mathfrak{S}^{p,q}_{{\cal D},\hat u,\hat w}
\end{align}
(it can be proved similarly as Lemma 3.4 in \cite{vas_har_tree}).

We have
\begin{align}
\label{wlq} \|\hat w\|_{l_q({\cal D}_{\hat \xi_m})}
=\|w\|_{l_q({\cal A}_{\hat\xi_m})}.
\end{align}
Hence,
$$
\frac{\|\hat w\|_{l_q({\cal D}_{\xi'})}}{\|\hat w\|_{l_q({\cal
D}_\xi)}} \stackrel{(\ref{wwsig})}{<}\sigma, \quad \xi \in {\bf
V}({\cal D}), \quad \xi' \in {\bf V}_1^{{\cal D}}(\xi).
$$
Applying Lemma \ref{sigma}, we get that
$$
\mathfrak{S}^{p,q}_{{\cal D}, \hat u,\hat w}
\underset{p,q}{\asymp} \sup _{1\le m\le m_*} \hat u(\hat \xi_m)
\|\hat w\|_{l_q({\cal D}_{\hat \xi_m})}
\stackrel{(\ref{hatuxim}),(\ref{wlq})}{=}\sup _{1\le m\le m_*}
\|u\|_{l_{p'}({\cal A}_{\hat \xi_m,\sigma})}\|w\|_{l_q({\cal
A}_{\hat\xi_m})}\le
$$
$$
\le \sup _{\xi \in {\bf V}({\cal A})}\|u\|_{l_{p'} ({\cal A}_{
\xi, \sigma})}\|w\|_{l_q({\cal A}_{\xi})}
\stackrel{(\ref{uw})}{\underset{p,q}{\lesssim}} \sup _{\eta \in
{\bf V}({\cal A})} \|u\|_{l_{p'}([\xi_0, \, \eta])}
\|w\|_{l_q({\cal A}_{\eta})}.
$$
This together with (\ref{spq}) completes the proof.
\end{proof}

\renewcommand{\proofname}{\bf Proof}

\section{Proof of Theorem \ref{lpq}}

The following result was proved by G. Bennett \cite{bennett_g}.

\begin{trma}
\label{hardy_diskr} {\rm (see \cite{bennett_g}).}  Let $1<p \le q<
\infty$, and let $\hat u=\{\hat u_n\}_{n\in \Z_+}$, $\hat w=\{\hat
w_n\}_{n\in \Z_+}$ be non-negative sequences such that
$$
M_{\hat u,\hat w}:=\sup _{m\in \Z_+}\Bigl(\sum \limits
_{n=m}^\infty \hat w_n^q\Bigr)^{\frac 1q}\Bigl( \sum \limits
_{n=0}^m \hat u_n^{p'}\Bigr)^{\frac{1}{p'}}<\infty.
$$
Let $\mathfrak{S}^{p,q}_{\hat u,\hat w}$ be the minimal constant
$C$ in the inequality
$$
\left(\sum \limits _{n=0}^\infty\left|\hat w_n\sum \limits
_{k=0}^n \hat u_kf_k\right|^q\right)^{1/q}\le C\left( \sum \limits
_{n\in \Z_+}|f_n|^p\right)^{1/p}, \quad \{f_n\}_{n\in\Z_+}\in l_p.
$$
Then $\mathfrak{S}^{p,q}_{\hat u,\hat w} \underset{p,q}{\asymp}
M_{\hat u,\hat w}$.
\end{trma}

\renewcommand{\proofname}{\bf Proof of Theorem \ref{lpq}}
\begin{proof}
By (\ref{spqh}), (\ref{ujwj}), (\ref{cvj}) and Theorem
\ref{sum_oper}, it is sufficient to prove that
$$
\hat{\mathfrak{S}}^{p,q}_{{\cal A},u,w} \underset{p,q, C_*,R,
R_0}{\lesssim} \sup _{j\in \Z_+} u_j\left(\sum \limits _{i\ge j}
w_i^q \frac{S(i)}{S(j)}\right)^{\frac 1q}.
$$

Denote by $\{\eta_{j,i}\}_{i\in I_j}$ the set of vertices ${\bf
V}_j^{{\cal A}}(\xi_0)$.

Let $h:{\bf V}({\cal A}) \rightarrow \R_+$, $\|h\|_{l_p({\cal
A})}=1$, $f(\xi)=\sum \limits _{\eta\le \xi} u(\eta)h(\eta)$. We
estimate from above the magnitude
$$\sum \limits _{\xi \in {\bf V}({\cal A})} w^q(\xi) \left(\sum
\limits _{\eta \le \xi} u(\eta) h(\eta)\right)^q = \sum \limits
_{\xi \in {\bf V}({\cal A})} w^q(\xi)f^q(\xi).$$

Let $n\in \Z_+$. Denote by $X_n$ the disjoint union of intervals
$\Delta_{n,i}$ ($i\in I_n$) of unit length. Let $\mu_n$ be a
measure on $X_n$ such that $\mu_n(\Delta_{n,i})=1$ and the
restriction of $\mu_n$ on $\Delta_{n,i}$ is the Lebesgue measure.
Then
\begin{align}
\label{munxn} \mu_n(X_n)={\rm card}\, {\bf V}_n^{\cal A}(\xi_0)
\stackrel{(\ref{cvj}), (\ref{s01})} {\underset{C_*}{\asymp}} S(n).
\end{align}

We define the function $\varphi:X_n \rightarrow \R$ by
\begin{align}
\label{vrphdni} \varphi|_{\Delta_{n,i}} = f(\eta_{n,i}).
\end{align}

Given $0\le k\le n$, $s\in I_k$, we set $Q_{k,s}=\cup
_{\eta_{n,i}\ge \eta_{k,s}}\Delta_{n,i}$. Then
\begin{align}
\label{muqks} \mu_n(Q_{k,s}) ={\rm card}\, \{i\in I_n:\;
\eta_{n,i}\ge \eta_{k,s}\}
\stackrel{(\ref{cvj})}{\underset{C_*}{\asymp}} \frac{S(n)}{S(k)}.
\end{align}
Let
\begin{align}
\label{pkek} P_k\varphi|_{Q_{k,s}} =\sum \limits _{\eta<
\eta_{k,s}} u(\eta)h(\eta), \quad E_k\varphi=\varphi-P_k\varphi.
\end{align}
In addition, we set
\begin{align}
\label{en1} P_{n+1}\varphi:=\varphi, \quad E_{n+1}\varphi:=0.
\end{align}

Let us estimate from above the value
\begin{align}
\label{sliin} \sum \limits _{i\in I_n} |f(\eta_{n,i})|^q
\stackrel{(\ref{vrphdni})}{=} \|\varphi\|^q_{L_q(X_n)}.
\end{align}
To this end, we argue similarly as in \cite[Theorem
4.1]{i_irodova1}, \cite{i_irodova}.

Denote by $\varphi^*:[0, \, \mu_n(X_n)] \rightarrow \R_+$ the
non-increasing rearrangement of the function $|\varphi|$. Then for
any $t\in [0, \, \mu_n(X_n)]$ there exists a set $Y_t\subset X_n$
such that
\begin{align}
\label{yt} \mu_n(Y_t)=t, \quad \forall x\in Y_t \;\;\;\;
|\varphi(x)|\ge \varphi^*(t).
\end{align}
Then
\begin{align}
\label{varphilp} \varphi^*(t)\le \frac{\|\varphi\|_{L_p(Y_t)}}
{[\mu_n(Y_t)]^{1/p}}.
\end{align}

Let $k\in \Z_+$,
\begin{align}
\label{tsn} \frac{\mu_n(X_n)}{S(k+1)} < t\le
\frac{\mu_n(X_n)}{S(k)}.
\end{align}
We set $k_n=\min\{k, \, n+1\}$. Then
$$\|\varphi\|_{L_p(Y_t)} \stackrel{(\ref{pkek}),(\ref{en1})}{\le} \|E_{k_n}\varphi\|_{L_p(Y_t)} +\sup
_{x\in X_n} |P_{k_n}\varphi(x)| [\mu_n(Y_t)]^{1/p}.$$ This
together with (\ref{varphilp}) yields that
$$
\varphi^*(t)\le \|E_{k_n}\varphi\| _{L_p(Y_t)}[\mu_n(Y_t)]
^{-1/p}+\sup _{x\in X_n} |P_{k_n}\varphi(x)|.
$$
From (\ref{sj}), (\ref{yt}) and (\ref{tsn}) we obtain
\begin{align}
\label{fst} \varphi^*(t)\le \|E_{k_n}\varphi\| _{L_p(Y_t)}
[\mu_n(X_n)]^{-1/p} [S(k)]^{1/p} +\sup _{x\in X_n}
|P_{k_n}\varphi(x)|.
\end{align}

Let $k_n\ge 1$ (otherwise, $P_{k_n}\varphi=0$), and let $x\in
Q_{k_n-1,s}$. For any $0\le j\le k_n-1$ there exists the unique
$l_j\in I_j$ such that $\eta_{j,l_j}\le \eta_{k_n-1,s}$. Then
\begin{align}
\label{pj1} (P_{j+1}\varphi -P_j\varphi) |_{Q_{j,l_j}}
\stackrel{(\ref{pkek}), (\ref{en1})}{=} {\rm const}, \quad 0\le
j\le k_n-1.
\end{align}
Observe that $P_0\varphi=0$. Hence,
$$
|P_{k_n}\varphi(x)| \le \sum \limits _{j=0}^{k_n-1}
|P_{j+1}\varphi(x) -P_j\varphi(x)| \le \sum \limits _{j=0}^{k_n-1}
\|P_{j+1}\varphi -P_j\varphi\|_{C(Q_{j,l_j})}
\stackrel{(\ref{pj1})}{=}
$$
$$
= \sum \limits _{j=0}^{k_n-1} \|P_{j+1}\varphi -P_j
\varphi\|_{L_p(Q_{j,l_j})} [\mu_n(Q_{j,l_j})]^{-1/p}
\stackrel{(\ref{muqks}), (\ref{pkek}), (\ref{en1})}{\underset{p,
C_*}{\lesssim}}
$$
$$
\lesssim\sum \limits _{j=0}^{k_n-1}
\left[\frac{S(j)}{S(n)}\right]^{1/p} (\|E_j \varphi\|
_{L_p(Q_{j,l_j})} +\|E_{j+1} \varphi\| _{L_p(Q_{j,l_j})})
\stackrel{(\ref{sj})}{\underset{p, R,R_0}{\lesssim}}
$$
$$
\lesssim \sum \limits _{j=0}^{k_n}
\left[\frac{S(j)}{S(n)}\right]^{1/p} \|E_j \varphi\| _{L_p(X_n)}
\stackrel{(\ref{en1})}{=} \sum \limits _{j=0}^{\min\{k,n\}}
\left[\frac{S(j)}{S(n)}\right]^{1/p} \|E_j \varphi\| _{L_p(X_n)}.
$$
This together with (\ref{munxn}) and (\ref{fst}) implies that
\begin{align}
\label{vrphstt} \varphi^*(t) \underset{p,C_*,R,R_0}{\lesssim} \sum
\limits _{j=0}^{\min \{k, \, n\}}
\left[\frac{S(j)}{S(n)}\right]^{1/p} \|E_j \varphi\| _{L_p(X_n)}
\quad \text{if}\;\;\frac{\mu_n(X_n)}{S(k+1)} < t\le
\frac{\mu_n(X_n)}{S(k)}.
\end{align}
Therefore,
\begin{align}
\label{intllxn} \begin{array}{c} \int \limits
_{X_n}|\varphi(x)|^q\, dx =\int \limits _0^{\mu_n(X_n)}
|\varphi^*(t)|^q\, dt \stackrel{(\ref{s01}), (\ref{sj}),
(\ref{munxn}), (\ref{vrphstt})}{\underset{p,q,
C_*,R,R_0}{\lesssim}}
\\
\lesssim \sum \limits _{k=0}^n \frac{S(n)}{S(k)} \left(\sum
\limits _{j=0}^k \left[\frac{S(j)}{S(n)}\right]^{1/p} \|E_j
\varphi\| _{L_p(X_n)}\right)^q=:A.
\end{array}
\end{align}
We claim that
\begin{align}
\label{aest} A\underset{p,q, C_*,R,R_0}{\lesssim} \sum \limits
_{j=0}^n \left(\frac{S(j)}{S(n)}\right)^{\frac qp-1} \|E_j
\varphi\|^q_{L_p(X_n)}.
\end{align}
Indeed, let $\psi_j =\left(\frac{S(j)}{S(n)}\right)^{\frac
1p-\frac 1q} \|E_j\varphi\| _{L_p(X_n)}$. Then (\ref{aest})
follows from Theorem \ref{hardy_diskr} and the estimate
$$
\sup _{0\le k\le n} \left(\sum \limits
_{j=0}^k\left[\frac{S(j)}{S(n)}\right]^{q'/q}
\right)^{\frac{1}{q'}} \left(\sum \limits _{j=k}^n
\frac{S(n)}{S(j)} \right)^{\frac 1q}
\stackrel{(\ref{sj})}{\underset{p,q,R,R_0}{\lesssim}}
$$
$$
\lesssim \sup _{0\le k\le n} \left(\frac{S(k)}{S(n)}\right)^{\frac
1q} \cdot\left(\frac{S(n)}{S(k)}\right)^{\frac 1q}=1.
$$

We have
$$
\|E_j\varphi\|_{L_p(X_n)} \stackrel{(\ref{pkek}), (\ref{en1})}{=}
\left\| \sum \limits
_{i=j+1}^{n+1}(P_i\varphi-P_{i-1}\varphi)\right\|_{L_p(X_n)} \le
\sum \limits _{i=j+1}^{n+1}
\|P_i\varphi-P_{i-1}\varphi\|_{L_p(X_n)}.
$$
From (\ref{intllxn}) and (\ref{aest}) we get that
\begin{align}
\label{flqest} \begin{array}{c}
\|\varphi\|_{L_q(X_n)}^q\underset{p,q,R,R_0,C_*}{\lesssim} \sum
\limits _{j=0}^n \left(\frac{S(j)}{S(n)}\right)^{\frac qp-1}
\left(\sum \limits _{i=j+1}^{n+1}
\|P_i\varphi-P_{i-1}\varphi\|_{L_p(X_n)}\right)^q
{\underset{p,q,R,R_0}{\lesssim}} \\ \lesssim \sum \limits
_{j=1}^{n+1} \left(\frac{S(j)}{S(n)}\right)^{\frac
qp-1}\|P_j\varphi-P_{j-1}\varphi\|^q_{L_p(X_n)}.
\end{array}
\end{align}
The last inequality follows from Theorem \ref{hardy_diskr} and
(\ref{sj}) since
$$
\sup _{1\le j\le n+1} \left(\sum \limits_{i=1}^j
\left(\frac{S(i)}{S(n)}\right)^{\frac qp-1}\right)^{1/q}
\left(\sum \limits_{i=j}^{n+1}
\left(\frac{S(n)}{S(i)}\right)^{q'\left(\frac 1p-\frac
1q\right)}\right)^{1/q'} \stackrel{(\ref{sj})}
{\underset{p,q,R,R_0}{\lesssim}} 1.
$$

Further,
$$
\|P_j\varphi-P_{j-1}\varphi\|_{L_p(X_n)}^p
\stackrel{(\ref{pkek})}{=}\sum \limits _{s\in I_{j-1}} \int
\limits _{Q_{j-1,s}} \left|\sum \limits _{\eta \le \eta _{j-1,s}}
u(\eta)h(\eta)-\sum \limits _{\eta < \eta _{j-1,s}}
u(\eta)h(\eta)\right|^p\, d\mu_n \stackrel{(\ref{ujwj}),
(\ref{muqks})}{\underset{p,q,R,R_0}{\lesssim}}
$$
$$
\lesssim\sum \limits _{s\in I_{j-1}} u_{j-1}^p|h(\eta_{j-1,s})|^p
\frac{S(n)}{S(j-1)} =\frac{S(n)}{S(j-1)} u_{j-1}^p \sum \limits
_{s\in I_{j-1}} |h(\eta_{j-1,s})|^p.
$$
This together with (\ref{sj}) and (\ref{flqest}) implies that
$$
\|\varphi\|_{L_q(X_n)}^q\underset{C_*,p,q,R,R_0}{\lesssim}\sum
\limits _{k=0}^{n} \frac{S(n)}{S(k)} u_k^q \left(\sum \limits
_{s\in I_k} |h(\eta_{k,s})|^p\right)^{\frac qp}.
$$
Hence,
$$
\sum \limits _{n=0}^\infty w^q_n \sum \limits _{i\in I_n}
|f(\eta_{n,i})|^q
\stackrel{(\ref{sliin})}{\underset{p,q,R,R_0}{\lesssim}} \sum
\limits _{n=0}^\infty w^q_n\sum \limits _{k=0}^{n}
\frac{S(n)}{S(k)} u_k^q \left(\sum \limits _{s\in I_k}
|h(\eta_{k,s})|^p\right)^{\frac qp} =
$$
$$
=\sum \limits _{k=0}^\infty \left(\sum \limits _{s\in I_k}
|h(\eta_{k,s})|^p\right)^{\frac qp}u_k^q \sum \limits
_{n=k}^\infty w_n^q\frac{S(n)}{S(k)}\le
$$
$$
\le \left[\sup _{0\le k<\infty} u_k^q\sum \limits _{n=k}^\infty
w_n^q\frac{S(n)}{S(k)}\right] \sum \limits_{k=0}^\infty \left(\sum
\limits _{s\in I_k} |h(\eta_{k,s})|^p\right)^{\frac qp}.
$$
This completes the proof.
\end{proof}

\renewcommand{\proofname}{\bf Proof}

\end{document}